\documentclass[11pt]{article}

\usepackage{epsfig}
\usepackage{graphicx}
\usepackage{amsmath}
\usepackage{amssymb}
\usepackage{amsfonts}

\newtheorem{thm}{Theorem}
\newtheorem{prop}[thm]{Proposition}
\newtheorem{cor}[thm]{Corollary}
\newtheorem{lemma}[thm]{Lemma}

\newenvironment{pf}{\noindent\textbf{Proof}}{\par}

\def\qed{\hspace{5truemm}\mbox{$\Box$}}

\newcommand{\subgp}[1]{\langle{#1}\rangle}

\def\Title#1{\begin{center} {\Large \textbf{#1} } \end{center}}

\begin{document}

\Title{TOPOLOGY OF CAYLEY GRAPHS APPLIED TO INVERSE ADDITIVE PROBLEMS\\Dedicated to M. A. FIOL}
\medskip

\begin{raggedright}
\emph{
Yahya Ould Hamidoune\\                 %%Name of the first autor
Universit\'e Pierre et Marie Curie, E. Combinatoire\\             %%Institution
Paris}                             %%City
\end{raggedright}

\begin{abstract}
We present proofs of the basic isopermetric structure theory, obtaining some new simplified proofs.
As an application, we obtain simple descriptions for subsets $S$ of an abelian group with
$|kS|\le k|S|-k+1$  or $|kS-rS|- (k+r)|S|,$ where $1\le r \le k.$ These results may be applied to several questions
in Combinatorics and Additive Combinatorics (Frobenius Problem, Waring's problem in finite fields and Cayley graphs with a big diameter, ....).
\end{abstract}
\section{ Introduction}
The connectivity of a graph is just the smallest number of vertices disconnecting the graph. In order investigate
more sophisticated properties of graphs, several authors proposed generalizations of connectivity. The reader may
find details on this investigation in the chapter \cite{FabFio04}. Investigating the isoperimetric connectivity in Cayley
graphs is just one of the many facets of Additive Combinatorics. It is also one of the many facets of Network topology.
For space limitation, we concentrate on  Additive Combinatorics, but the reader may find details and a bibliography in 
the recent paper
\cite{HamLlaLop} concerning the other aspect. 

Let $\Gamma=(V,E)$ be a reflexive graph. The minimum of the objective function $|\Gamma (X)|-|X|,$ restricted to subsets
$X$ with $|X|\ge k$ and $|V \setminus \Gamma (X)|\ge k,$ the $k$-isoperimetric connectivity. Subsets  achieving the above minimum are called $k$-fragments. $k$-fragments with smallest cardinality are called $k$-atoms. It was proved by the author in \cite{Ham96}, that distinct $k$-atoms of $\Gamma$
 intersect in at most $k-1$ elements, if the size of the $k$-atom
of $\Gamma$ is not greater than the size of the $k$-atom
of $\Gamma^{-1}.$ Let $1\in S$ be a finite generating subset of a group $G$ such that the cardinality
$1$-atom of the Cayley graph defined by $S$ is not greater than the cardinality
$1$-atom of the Cayley graph defined by $S^{-1}.$ Then a $1$-atom $H$ containing $1$ is a subgroup.
The last  result applied to a group with a prime order is just the Cauchy-Davenport Theorem.
It has several other implications and leads to few lines proof for result having very tedious proof using the classical transformations. In particular, it was applied recently by the author \cite{hkft} to a problem of Tao \cite{t2}.

In the abelian case, things are much easier. Assume that $G$ is abelian and let $1\in H$ be a $k$-atom of the Cayley graph defined by $S$. If $k=1,$ then $H$ a subgroup (the condition involving $S^{-1}$ is automatically verified). In particular, there is a subgroup which is a $1$-fragment.
A maximal such a group is called an hyper-atom. Assuming now that $k=2$ and that $\kappa _2\le |S|-1.$
 It was proved in \cite{Ham97} that either $|H|=2$ or $H$
is a subgroup. It was proved also in \cite{Ham97} that either $S$ is an arithmetic 
progression or there is a non-zero subgroup which is a $1$-fragment,  if $|S|\le (|G|+1)/2$. Let $Q$ be a hyper-atom of $S$ and let  $\phi :G\mapsto G/Q$ denotes
the canonical morphism. The author proved in \cite{hkemp}  that
$\phi (S)$ is either an arithmetic progression or satisfies the sharp Vosper property (to be defined later) if $|S|\le (|G|+1)/2$.

Let $G$ be an abelian group and let $A,B$ be finite non-empty subsets of $G,$ with $|A+B|= |A|+|B|-1-\mu.$
Kneser's Theorem states that $\pi (A+B)\neq \{0\},$ where $\pi (A+B)=\{x: x+A+B=A+B\}.$ The hard Kemperman Theorem,
which needs around half a page to be formulated, describes recursively the subsets $A$ and $B$ if $\mu =1.$ Its classical proof requires around 30 pages. It was applied by Lev \cite{levkemp} to propose a dual description, that looks
easier to implement than Kemperman's description.

The above structure isoperimetric results were used in \cite{hkemp,hkemp+1} to explain the topological nature of Kemperman Theory and to give a shorter proof of it. Our method involve few technical steps and use some duality arguments and the strong isoperimetric property. We suspect that it
could be drastically simplified. In this paper, we shall verify this hypothesis for Minkowski sums of the form $rS-sS,$ obtaining very simple proofs and tight descriptions. This case covers almost all the known applications. Also, Modern Additive Combinatorics deals almost exclusively with $rS-sS$, c.f. \cite{tv}.

The organization of the paper is the following:

Section 2 presents the isoperimetric  tools, with complete proofs in order to make the paper self-contained. In particular, this section contains a proof of the fundamental property of $k$-atoms. 
 In Section 3, we start by showing the structure of $1$-atoms of arbitrary Cayley graphs. We then restrict ourselves to the abelian case. We give in this section   an new simplified proof for the structure theorem of $2$-atoms. We deduce from it 
the structure of hyper-atoms. In Section 4, we give easy properties of the decomposition modulo a subgroup which is a fragment.
Easy proofs of the Kneser's theorem and a Kemperman type result for $kA$ are then presented.

In Section 5, we investigate universal periods for $kS$ introducing a new object: the sub-atom.
It follows from a result by Balandraud \cite{balart} that $|TS|\le |T|+|S|-2$ implies that $T+S+K=T+S,$ where $K$ is the final kernel of $S$ (a subgroup contained in the atom of $S$ described in \cite{balart}). We shall prove that  the  $kS+M=kS,$ if $|kS|\le k|S|-k,$ where $M$ is the sub-atom. Clearly $K\subset M.$ The case $rS-sS,$ where $r\ge s\ge 1$, is solved easily in Section 5, by showing that one of the following holds\begin{itemize}
  \item $S$ is an arithmetic progression,
  \item $|sS-rS|\ge \min (|G|-1,(r+s)|S|),$
  \item $|H|\ge 2$ and $sS-rS+H=sS-rS,$ where $H$ is an  hyper-atom  of $S.$ 
  \end{itemize}

Readers familiar with Kemperman Theory could appreciate the simplicity of this result.
 In Section 6, we obtain the following description: 
 
 Let $k\ge 3$ be  an integer and let $0\in S$ be a finite generating  subset of an abelian group $G$ such that $S$ is not an arithmetic progression, $kS$ is aperiodic and $|kS|= k|S|-k+1.$
 Let $H$ be a hyper-atom of $S$ and let $S_0$ be a smaller $H$-component of $S.$
Then $(S\setminus S_0)+H=(S\setminus S_0)$  and  $|kS_0|=k|S_0|-k+1.$
Moreover  $\phi (S)$ is an arithmetic progression, where $\phi :G\mapsto G/H$ denotes
the canonical morphism.
  
Necessarily $|H|\ge 2,$ since $S$ is not an arithmetic progression.
\section{Basic notions}
 Recall a well known fact:

\begin{lemma}(folklore)
Let $a,b$ be elements of a group $G$ and let $H$ be a finite subgroup of $G.$
Let
 $A,B$  be  subsets of $G$
 such that  $A\subset aH$ and $B\subset Hb.$ If  $|A|+|B|>|H|,$
 then $AB=aHb$.
\label{prehistorical}
 \end{lemma}

Let  $H$ be a subgroup of an abelian group $G$. Recall that a {\em  $H$--coset} is a set of the $a+H$ for some $a\in G$.
The family $\{a+H; a\in G\}$ induces a partition of $G$. The trace of this partition on a subset $A$
will be called an {\em  $H$--decomposition}
of $A$.

By a graph, we shall mean a directed graph,
identified  with its underlying relation.  Undirected graphs are identified with symmetric graphs. We recall the definitions in this
context.

An ordered pair $\Gamma =(V,E),$ where $V$ is a set and  $E\subset V\times V,$ will be called a {\em graph}
or a {\em relation} on  $V.$ Let $\Gamma =(V,E)$ be a  graph and let $X\subset V.$ The {\em reverse} graph of $\Gamma$ is
the graph $\Gamma ^- =(V,E^-),$ where $E^-=\{(x,y): (y,x)\in E\}.$
The {\em degree} (called also outdegree) of a vertex $x$ is $$ d(x)=|\Gamma (x)|.$$
%Graphs could be infinite, but all graphs considered here have finite degrees.
The graph $\Gamma$ will be called {\em locally-finite}  if for all $x\in V,$ $|\Gamma (x)|$ and $|\Gamma ^-(x)|$ are finite.
The graph $\Gamma$ is said to be {\em $r$-regular} if $|\Gamma (x)|=r,$ for every $x\in V.$
The graph $\Gamma$ is said to be {\em $r$-reverse-regular} if $|\Gamma ^-(x)|=r,$ for every $x\in V.$
The graph $\Gamma$ is said to be {\em $r$-bi-regular} if it is {\em $r$-regular} and {\em $r$-reverse-regular}.

 \begin{itemize}
  \item The minimal degree of $\Gamma$ is defined as $\delta ({\Gamma})= \min \{ |\Gamma (x)| : x\in V\}.$
  \item We write $\delta _{\Gamma ^-}= \delta _{-} (\Gamma).$
  \item The boundary of  $X$  is defined as $\partial _{\Gamma}(X)= \Gamma (X)\setminus X.$
  \item The exterior of  $X$  is defined as
 $\nabla _{\Gamma}(X)= V\setminus \Gamma (X).$
  \item We shall write $\partial ^{-}_{\Gamma } =\partial _{\Gamma ^{-}}.$ This subset will be called the
  {\em reverse-boundary} of  $X.$
  \item  We shall write $\nabla ^{-}_{\Gamma } =\nabla _{\Gamma ^{-}}.$
\end{itemize}

In our approach, $\Gamma (v)$ is just the image of $v$ by the relation $\Gamma$ and  $\Gamma^-(v)$ requires no definition
since $\Gamma ^-$ is defined in Set Theory as the reverse of $\Gamma$.

An automorphism  of a graph $\Gamma =(V,E)$ is a permutation $f$ of $V$ such that $f(\Gamma (v))=\Gamma (f(v)),$
for any vertex $v.$
A graph $\Gamma =(V,E)$ is said to be {\em vertex-transitive} if for any ordered pair of vertices there is an automorphism
mapping the first one to the second.

Let  $ A,B$ be
subsets of a group $ G $. The {\em Minkowski product} of $ A$ with $B$ is defined as
$$AB=\{xy \ : \ x\in A\  \mbox{and}\ y\in
  B\}.$$

Let  $S$  be a
subset of   $G$. The subgroup
generated by $S$ will be denoted by $\subgp{S}$.
The graph $(G,E),$ where  $ E=\{ (x,y) : x^{-1}y \
\in S \}$ is called a {\it Cayley graph}.  It will  be denoted by
$\mbox{Cay} (G,S)$.
Put $\Gamma =\mbox{Cay} (G,S)$   and  let   $F \subset G $.
Clearly
 $\Gamma (F)=FS .$
One may check easily that left-translations are automorphisms of Cayley graphs.
In particular,  Cayley  graphs are   vertex-transitive.

Let $\Gamma=(V,E)$ be a reflexive graph.
%In other words $\Gamma$ is reflexive relation or a directed graph with loops at every vertex.
 We shall investigate the {\em boundary operator}  $\partial _\Gamma : 2^V \rightarrow 2^V.$
When the context is clear, the reference to $\Gamma$ will be omitted.
Since $\Gamma$ is reflexive, we have in this case $ |\partial (X)| =|\Gamma (X)|-| X|.$

Let  ${\mathcal A} \subset 2^V$ be a family of finite subsets of $V$.
 We define the {\em connectivity} of ${\mathcal A}$ as  $$\kappa ({\mathcal A}) =
\min  \{|\partial (X)|\   :  \ \ X\in {\mathcal A}\}.$$

An $X\in  {\mathcal A}$ with $\kappa ({\mathcal A}) =
|\partial  (X)|$ will be called a {\em fragment}.

A fragment with a minimal cardinality will be called an {\em atom}.

%There are interesting families. But we consider here only the family of $k$-separable sets defined below:

Put
$${\mathcal S}_k(\Gamma)=\{ X :  k\le |X|<\infty \ \text{and} \ |\Gamma (X)|\le |V|-k \ \}.$$
%The elements of ${\mathcal S}_k$ are called {\em $k$-separable} subsets.

We shall say that $\Gamma$ is {\em $k$-separable} if ${\mathcal S}_k(\Gamma)\neq \emptyset.$ In this case, we write $$\kappa _k(\Gamma)=\kappa ({\mathcal S}_k) .$$

By a $k$-{\em fragment} (resp.  $k$-{\em atom}), we shall mean a {fragment} (resp.  { atom}) of ${\mathcal S}_k$.
A $k$-fragment of $\Gamma^{-1}$ is sometimes called a  $k$-{\em negative}) fragment. This notion was introduced by the author in \cite{Ham96}. A relation $\Gamma$ will be called $k$-{\em faithful} if $|A|\le |V\setminus \Gamma (A)|,$
 where $A$ is a $k$-atom of $\Gamma$.
 By a {\em fragment} (resp.  {\em atom}), we shall mean a $1$-{fragment} (resp.  $1$-{ atom}).

The following lemma is immediate from the definitions:

\begin{lemma} \cite{Ham96}\label{degenerate}{ Let $k\ge 2$ be an integer. A  reflexive
  locally finite   $k$-separable graph  $\Gamma =(V,E)$ is a $k-1$-separable graph, and moreover
  $\kappa _{k-1}\le \kappa _k.$

  }
  \end{lemma}

Recall the following easy fact:
\begin{lemma}\label{katomdegree}\cite{Ham96}
Let $\Gamma =(V,E)$ be a  locally-finite $k$-separable graph and let $A$ be a $k$-atom with $|A|>k.$
 Then $\Gamma ^-(x)\cap A\neq \{x\},$ for every $x\in A.$ %In particular,  $A$ contains a directed cycle.
\end{lemma}

\begin{pf} We can not have $\Gamma ^-(x)\cap A= \{x\},$ otherwise $ A\setminus \{x\}$
would be a $k$-fragment.
\qed\end{pf}

The next lemma contains useful duality relations:

\begin{lemma} \cite{Ham99}\label{negative}{Let  $X$ and $Y$ be $k$-fragments of a reflexive
  locally finite   $k$-separable graph $\Gamma =(V,E)$.  Then
 \begin{align}
\partial ^{-}  (\nabla (X))&=
\partial (X),\label{eqduall}\\
\nabla^- (\nabla (X))&={X}, \label{eqdualf}
 \end{align}

}\end{lemma}
\begin{pf} Clearly, $ \partial (X) \subset \partial ^{-} (\nabla (X))$

We must have $ \partial (X) = \partial ^{-} (\nabla (X))$, since
otherwise there is a $y\in \partial ^{-} (\nabla (X))
\setminus \partial (X).$ It follows that
$|\partial (X\cup \{y\})|\le |\partial (X)|-1$, contradicting the
definition of $\kappa _k.$ This proves (\ref{eqduall}).

Thus $\Gamma ^{-}
(\nabla (X))=\nabla (X)\cup
\partial ^{-}(\nabla (X)) =\nabla (X)\cup
\partial (X)=V\setminus X,$
and hence (\ref{eqdualf}) holds.

Let $\Gamma =(V,E)$ be a  reflexive graph. We shall say that
$\Gamma$ is a {\em Cauchy graph}  if $\Gamma$ is  non-$1$-separable or if $\kappa _1(\Gamma)=\delta-1.$ h  We shall say that
$\Gamma$ is a {\em reverse-Cauchy graph} if $\Gamma ^-$ is a Cauchy graph.

Clearly, $\Gamma$ is a Cauchy graph if and only if for every
$X\subset V$ with $|X|\ge 1$,
 $$|\Gamma(X)|\ge \min \Big(|V|, |X|+\delta-1\Big).$$
 \qed\end{pf}

 \begin{lemma} \cite{Ham96}\label{finiteg}{Let $\Gamma =(V,E)$ be a reflexive
finite   $k$-separable graph and let $X$  be a subset of $V.$
Then \begin{equation}
\kappa _k= \kappa _{-k}.\label{kk-}\end{equation}  Moreover,
\begin{itemize}
  \item[(i)]   $X$ {is a } $k$-{fragment } {if and only if } $\nabla(X)$ {is a } $k$-{reverse-fragment,}
  \item[(ii)]  $\Gamma$ is a Cauchy graph if and only if it is a reverse-Cauchy graph.

\end{itemize}
}\end{lemma}
\begin{pf}
Observe that a finite graph is $k$-separable if and only if its reverse
is $k$-separable.
Take a
$k$-fragment $X$  of $\Gamma$. We have clearly $\partial _{-}
(\nabla (X))\subset
\partial (X)$.
Therefore $$\kappa _k(\Gamma )\ge |\partial (X)|\ge |\partial ^{-}
(\nabla (X))|\ge \kappa _{-k}.$$ The reverse inequality of (\ref{kk-}) follows similarly or by duality.

Suppose that $X$ is a $k$-fragment.
By (\ref{eqduall}) and (\ref{kk-}), $|\partial _{-}
(\nabla (X))|=
|\partial (X)|=\kappa _{k}=\kappa _{-k},$ and hence $\nabla (X)$ is a revere $k$-fragment. The other
implication of (i) follows easily.
Now (ii)  follows directly from the definitions.% and (\ref{kk-}).
\qed\end{pf}

\begin{thm}\cite{Ham96}

Let $\Gamma =(V,E)$ be a  reflexive  locally-finite $k$-faithful
$k$-separable graph. Then
the intersection of two distinct $k$-atoms $X$ and $Y$ has a cardinality less than  $k.$
Moreover, any locally-finite
$k$-separable graph is either  $k$-faithful or reverse   $k$-faithful.
\label{katom} \end{thm}
\begin{pf}
{\begin{center}
\begin{tabular}{|c||c|c|c|c|}
\hline
$\cap $&  $Y$ & $\partial (Y) $ &  $\nabla (Y)$ \\
\hline
\hline
$X$&  $R_{11}$ & $R_{12}$ & $R_{13}$ \\
\hline
$\partial (X) $& $R_{21}$& $R_{22}$ & $R_{23}$ \\
\hline
$\nabla (X)$ & $R_{31}$  & $R_{32}$ & $R_{33}$ \\
\hline
\end{tabular}

\end{center}}

Assume that $|X\cap Y|\ge k.$
 By the definition of $\kappa_k,$
\begin{align*}  |R_{21}|+|R_{22}|+ |R_{23}| &= \kappa_k
\\&\le |\partial (X\cap Y)|\\
&=  |R_{12}|+|R_{22}|+ |R_{21}|,
\end{align*}
and hence \begin{equation}|R_{23}|\le
|R_{12}|. \label{2312}\end{equation}

Thus,
 \begin{align*}|\nabla (X)\cap \nabla (Y)|&= |\nabla (Y)|-|R_{13}|- |R_{23}|\\
  &\ge |Y|-|R_{13}|- |R_{12}|\\&=|X|-|R_{13}|- |R_{12}|=|R_{11}|\ge k.
\end{align*}

Thus,

\begin{align*}  |R_{12}|+|R_{22}|+ |R_{32}| &=  \kappa_{k}
\\&\le |\partial (X\cup Y)|\\
&\le  |R_{22}|+|R_{23}|+ |R_{32}|,
\end{align*}
and hence $|R_{12}|\le |R_{23}|,$ showing that $|R_{12}|= |R_{23}|.$

 It follows that
   $$\kappa _k\le |\partial (X\cap Y)|\le  |R_{12}|+|R_{22}|+ |R_{21}|\le  |R_{12}|\le |R_{23}|+|R_{22}|+ |R_{21}|=\kappa _k,$$
showing that $X\cap Y$ is a $k$-fragment, a contradiction.

The fact that a locally-finite
$k$-separable graph is either  $k$-faithful or reverse   $k$-faithful follows by Lemma \ref{finiteg}.
 \qed\end{pf}

\section{A structure Theory for atoms}
In the sequel, we identify $\mbox{Cay} (\subgp{S},S)$ with $S,$ if $0\in S.$
We shall even work with subsets not containing $1$. By $\kappa _k(S)$ we shall mean
$\kappa _k(S-a)=\kappa _k(\mbox{Cay} (\subgp{S-},S-a)),$ for some $a\in S.$ As an exercise,
the reader may check that this notion does not depend on a particular choice of $a\in S.$
%The reader should have in mind that
\begin{thm}\label{Cay}\cite{hejc2}

Let $1\in S$ be a finite proper generating subset of a group $G$. Let $1\in H$ be a $1$-atom of $S.$
\begin{itemize}
  \item[(i)] If $S$ is $1$--faithful, then $H$ is a subgroup. Moreover  $|H|$ divides $\kappa _1(S).$
\item[(ii)] If $G$ is abelian and if $S$ is $k$-separable, then  $S$ is $k$-faithful.
  \item[(iii)]
 If $G$ is abelian, then  $H$ is a subgroup.
\end{itemize}
\end{thm}

\begin{pf} Take an element $x\in H.$ Clearly $xH$ is a $1$-atom. Since $(xH)\cap H\neq \emptyset,$ we have by Theorem \ref{katom}, $xH=H.$
Since $H$ is finite, $H$ is a subgroup. Now  $\kappa _1(S)=|HS|-|H|,$ showing the last part of (i).

If $G$ is abelian, then $\mbox{Cay} (G,S)$ is isomorphic to
$\mbox{Cay} (G,-S),$ and hence $S$ is $k$-faithful if $S$ is $k$-separable. Now (iii) follows by combining (i) and (ii).
  \qed\end{pf}
\begin{thm} { \cite{Ham97}\label{2atom} Let  $S$  be a finite generating
 $2$--separable subset of an abelian group $G$ with $0\in S$ and $\kappa _2 (S) \leq |S|-1$.
If
 $0\in H$ is a $2$--atom with $|H|\ge 3,$ then  $H$ is  a subgroup.
\label{2atomejc} }
\end{thm}

\begin{pf}
The proof is by induction. Assume first that $H+Q=H,$ where $Q$ is a non-zero subgroup.
For every, $x\in H,$ we have $|(H+x)\cap H|\ge |x+Q|\ge 2.$ By Theorems \ref{Cay} and  \ref{katom}, $H+x=H.$ It follows that
$H$ is a subgroup. Assume now that $H$ is aperiodic. Let us first show that $\kappa _1(H)=|H|-1.$ Suppose the contrary and take a
$1$-atom $L$ of $H$ with $0\in L.$ By Theorem \ref{Cay}, $L$ is a subgroup and $|L|\le \kappa _1(H).$ Take a nonzero element $y\in L.$ We have
$|H\cup (y+H)|\le |L+H|= |L|+\kappa_1 (H)\le 2\kappa _1(H)\le 2|H|-4.$ Thus, $|H\cap (y+H)|\ge 2,$ and hence
$y+H=H,$ by Theorem \ref{katom}.

%{\bf Case} 1. $H$ generates a proper subgroup $N$ of $G.$
Take an $N$-decomposition $S=\bigcup \limits_{1\le i\le s}
S_i,$ with $|S_1+H|\le \cdots \le |S_s+H|.$ Without loss of generality, we may take $0\in S_1.$ We have necessarily $s\ge 2.$
We must have $|S_i|=|N|,$ for all $i\ge 2.$ Suppose the contrary. By the definition of $\kappa _1,$ we have $|S_1+H|\ge |S_1|+\kappa _1 (H)=  |S_1|+|H|-1.$
We have also, since $H$ generates $N,$ $|S_i+H|\ge |S_1|+1.$
Thus, $|S+H|\ge   |S|+|H|-1+1\ge |S|+|H|,$ a contradiction.
Now we have $|X+S|=|S\setminus S_1|+|X+S_1,$ for any subset $X$ of $N.$
In particular, $H$ is a $2$-atom of $S_1.$ If $|S_1|<|S|,$ the result holds by Induction. It remains to consider the case $s=1.$

%Now we have $|N+S|-|N|=|S+H|-|S_1+H|\le |H+S|-|H|=\kappa _2(S)$ and hence $N=H,$ a contradiction.
%{\bf Case} 2. $H$ generates  $G.$
The relation $|H+S|-|H|\le |S|-1$ implies that $\kappa_2 (H)\le |H|-1.$
By Lemma \ref{katomdegree}, for every $x\in H,$ there $s_x\in S\setminus \{0\},$ with $x-s_x\in H.$ We must have $$|H|\le |S|-1,$$ otherwise there are distinct elements $x,y\in H$ and an element $s\in S\setminus \{0\}$ such that $x-s,y-s\in H.$ It follows that $|(H+s)\cap H|\ge 2.$ By Theorem \ref{katom}, $H+s=H,$ a contradiction.

Let $0\in M$ be a $2$-atom of $H.$ Take a non-zero element $a\in M.$
Since $\kappa _2(H)=|M+H|-|M|,$ $|M|$ divides $\kappa _2(H)$ if $M$ is a subgroup.
Thus, the Induction hypothesis implies that $|M|\le |H|-1.$ Since $|M+H|\le |M|+\kappa _2(H)\le 2|H|-2,$ we have $|H\cap (H+a)|\ge 2.$
  By Theorem \ref{katom}, $H+a=H,$ a contradiction.
\qed\end{pf}
%A short proof of this result is given in \cite{hiso2007}.

\begin{thm} { [\cite{Ham97},Theorem 4.6]
Let  $S$  be a $2$--separable finite
subset of an
 abelian group $G$ such that $0\in S$, $|S|\leq (|G|+1)/2$ and $\kappa _2 (S) \leq |S|-1$.

 If  $S$ is not an arithmetic progression, then there
 is a subgroup  which is a $2$--fragment of $S$.
\label{vosper0}                      }
\end{thm}

\begin{pf}

Suppose that $S$ is not an arithmetic progression.

 Let $H$ be a  $2$--atom such that $0\in H$. If $\kappa
_2\leq |S|-2$, then clearly $\kappa _2=\kappa _1$ and $H$
is also a $1$--atom. By Theorem \ref{Cay}, $H$ is a subgroup.
Then we may assume $$\kappa _2(S)=|S|-1.$$  By Theorem
\ref{2atomejc}, it would be enough to consider the case $|H|=2$, say
$H=\{0,x\}$. Put $N=\subgp{x}.$

Decompose $S=S_0\cup \cdots \cup S_j$ modulo $N$, where $|S_0+H|\le
|S_1+H| \le \cdots \le |S_j+H|.$ We have $|S|+1=|S+H|=\sum
\limits_{0\le i \le j}|S_i+\{0,x\}|.$

Then $|S_i|=|N|$, for all $i\ge 1$.  We have $j\ge 1$, since
otherwise $S$ would be an arithmetic progression. In particular,  $N$
is finite.
 We have
$|N+S|<|G|$, since otherwise   $|S|\ge |G|-|N|+1\ge
\frac{|G|+2}{2},$ a contradiction.

Now  \begin{eqnarray*} |N|+|S|-1&=&|N|+\kappa _2(S)\\&\le& |S+N|=
(j+1)|N|\\&\le& |S|+|N|-1, \end{eqnarray*}

 and hence $N$ is a $2$-fragment.
\qed\end{pf}

Theorem \ref{vosper0} was used to solve Lewin's
Conjecture on the  Frobenius number
 \cite{hacta}.

A $H$--decomposition  $A=\bigcup \limits_{i\in I}
A_i$  will be called
 a $H$--{\em modular-progression } if it is an arithmetic progression modulo $H$.

%In this section, we investigate  the  notion of a hyper-atom introduced in .
Recall that $S$ is a Vosper subset if and only if  $S$ is non
$2$--separable or if $\kappa _2(S)\ge |S|$.
%Theorem \ref{hyperatom} is one of the main results of this paper.

\begin{thm}\cite{hkemp}\label{hyperatom}
Let $S$ be a finite generating subset of an abelian group $G$ such
that $0 \in S,$  $| S | \leq (|G|+1)/2$ and $\kappa _2 (S)\le
|S|-1.$ Let $H$ be a hyper-atom of $S$. Then

$\phi (S)$ is either an arithmetic progression or a Vosper
subset, where $\phi$ is the canonical morphism from $G$ onto
$G/H$.
\end{thm}

\begin{pf}

Let us show that \begin{equation}\label{referee}
2|\phi (S)|-1\le \frac{|G|}{|H|}.\end{equation}
 Clearly we may
assume that $G$ is finite.

Observe that $2|S+H|-2|H|\le 2|S|-2< |G|.$ It follows, since$|S+H|$
is a multiple of $|H|$, that $2|S+H|\le  |G|+|H|,$ and hence (\ref{referee}) holds.

 Suppose now that  $\phi (S)$ is not a Vosper subset. By the
 definition of a Vosper subset,  $\phi (S)$ is $2$--separable and $\kappa _2(\phi(S))\le |\phi(S)|-1.$

Let us show that   $\phi(S)$
has no  $1$--fragment $M$ which is a non-zero subgroup.
Assuming the contrary. We have $|\phi(\phi ^{-1}(M)+S)|=|M+\phi (S)|\le |M|+|\phi (S)|-1.$
Thus,
$|\phi ^{-1}(M)+S|\le |\phi ^{-1}(M)|+|H|(|\phi (S)|-1)=|\phi ^{-1}(M)|+\kappa _1 (S).$ It follows that
$\phi ^{-1}(M)$ is a $1$-fragment. By the maximality of $H,$ we have $|M|=1,$ a contradiction.
 By (\ref{referee}) and
Theorem \ref{vosper0}, $\phi(S)$ is an arithmetic progression.\qed\end{pf}

%%%%%%%%%%%%%%%%%%%%%%%%%%%%%%%%%%%%%%%%%%%%

\section{Decomposition modulo a fragment}

Let  $H$ be a subgroup of an abelian group $G$. Recall that a {\em $H$--coset} is a set of the $a+H$ for some $a\in G$.
The family $\{a+H; a\in G\}$ induces a partition of $G$. A non-empty set of the form $A\cap (x+H)$
will be called a {\em  $H$-component}
of $A$. The partition of $A$ into its $H$-components
will be called a {\em  $H$-decomposition}
of $A.$ By a {\em smaller} component, we shall mean a component with a smallest cardinality.

Assume now that $H$ is $1$-fragment and take a $H$-decomposition $S=S_0\cup \dots \cup S_u,$ with
$|S_0|\le \dots \le |S_u|.$

We have $ |S|-1
\ge \kappa (S)=|H+S|-|H|.$

It follows that for $i\ge 1,$ we have $$2|H|-|S_0|-|S_i|\le |H+S|-|S|\le |H|-1,$$ and hence
$|S_0|+|S_i|\ge |H|+1.$ In particular,

for all $(i,j)\neq (0,0),$   $|S_i|+|S_j|\ge |H|+1,$  hence $$S_i+S_j+H=S_i+S_j,$$ by
Lemma \ref{prehistorical}.

Thus $$(S\setminus S_0)+S=(S\setminus S_0)+H+S.$$

Similarly $$((S\setminus S_0))-S=(S\setminus S_0)+H-S.$$
Since $S_0-S_0\subset S_1-S_1=H,$ we have $$S-S+H=S-S.$$

In particular,
$(kS\setminus kS_0)+H=kS\setminus kS_0.$

%It follows also that $\phi ((S\setminus H)+(k-1)S)\neq G/H.$
%Thus $u+|\phi (k-1)S)|\le |G|/|H|.$ Thus $|\phi ((k-1)S)|\ge (u+1)+(k-2)u=1+(k-1)u.$ In particular,
%$ 1+ku\le |G|/|H|.$ By iterating the isoperimetric inequality, and by Lemma, we have
%$|\phi (kS)|\ge \min (|G|/|H|, (u+1)+(k-1)u)=1+ku.$ %We may assume that  $|\phi (kS)| =1+ku,$ otherwise
%Now we have $|kS|\ge ku|H|+|kS_0|\ge ku|H|+k|S_0|-k+1\ge k|S|-k+1.$

%Assume first that $\kappa (S)=|S|-1.$
%Since $kS\neq G,$ we have $|kS|\ge |S|+(k-1)\kappa (S)=k|S|-k+1,$ and the result holds. So we may assume
%that  $ |S|-2\ge \kappa (S)=|H+S|-|H|.$

\begin{prop}
 Let $S_0$ denotes a smaller $H$-component of $S,$ where $H$ is a non-zero subgroup fragment.
We have $S-S+H=S-S.$
Let $2\le k$ be an integer. Then $(S\setminus S_0)+(k-1)S$ is $H$-periodic subset with cardinality at least
$ \min (|G|,k|S+H|-k|H|).$
If   $kS+H\neq kS,$ then $|S_1|>|H|/2\ge |S_0|,$ and $|kS|\ge k|S+H|-k|H|+|kS_0|.$
Moreover $kS_0$ is aperiodic if  $kS$ is aperiodic.

\label{decomp}
 \end{prop}

\begin{pf}
The first part was proved above. By the definition of $\kappa,$  we have  $|(S\setminus S_0)+(k-1)S|=|(S\setminus S_0)H+(k-1)S|\ge
u|H|+(k-1)\kappa =k|S+H|-k|H|.$

Assume now that $kS+H\neq kS.$  we have $kS_0\neq kS_0+H,$
and hence $2S_0\neq 2S_0+H,$ since $(S\setminus S_0)+(k-1)S$ is $H$-periodic.
By Lemma \ref{prehistorical}, $|H|/2\ge |S_0|.$
We have now $|S_1|\ge |H|+1-|S_0|\ge |H|/2+1.$ We must also have
$kS_0\cap  ((S\setminus S_0)+(k-1)S)=\emptyset.$ Thus,
$|kS|\ge |(S\setminus S_0)+(k-1)S|+|kS_0|\ge k|S+H|-k|H|+|kS_0|.$

Assume now that $kS$ is aperiodic. Since  $(S\setminus S_0)+(k-1)S$ is $H$-periodic and since  the period of $kS_0$ is a subgroup of $H$, necessarily $kS_0$ is aperiodic.
 \qed\end{pf}

\begin{cor} ( Kneser, \cite{knesrcomp})

Let $k$ be a non-negative integer and let $S$ be a finite  subset of an abelian group $G$. If $kS$ is aperiodic, then
$|kS|\ge k|S|-k+1$
\label{ejc2} \end{cor}

\begin{pf}
 Let $H$ be a $1$-atom containing $0.$ By Theorem \ref{Cay},  $H$ is subgroup. Let $S_0$ denotes a smaller $H$-component of $S.$ Without loss of generality we may assume that $0\in S_0.$  We may assume $\kappa (S)\le|S|-2
,$ since otherwise $|kS|\ge |S|+(k-1)\kappa (S)=k|S|-k+1,$ and the result holds.

By Proposition \ref{decomp}, $kS_0$ is aperiodic.
By the Induction hypothesis and Proposition \ref{decomp}, $|kS|=|kS_0|+ (k-1)(|S+H|-|H|)\ge k|S_0|-k+1+(k-1)(|S+H|-|H|)\ge k|S|-k+1.$
 \qed\end{pf}

We shall now complete Proposition \ref{decomp} in order to deal with the critical pair Theory.

\begin{prop}
Let $2\le k$ be an integer.
 Let $S_0$ denotes a smaller $H$-component of $S,$ where $H$ is a non-zero subgroup fragment $kS+H\neq kS.$ Assume moreover that $kS$ is aperiodic and  $|kS|= k|S|-k+1.$
Then

 \begin{itemize}
   \item[(i)]  $kS_0$ is aperiodic,
   \item [(ii)]   $|kS_0|= k|S_0|-k+1,$
   \item [(iii)]   $(S\setminus S_0)+H=S\setminus S_0$ and
   \item [(iv)]  $|k(S+H)|= k|S+H|-k|H|+|H|.$
 \end{itemize}

\label{crdecomp}
 \end{prop}

\begin{pf}
(i) follows by Proposition \ref{decomp}. By Kneser Theorem and Proposition \ref{decomp},
 \begin{align*}
k|S|-k+1=|kS|&\ge |kS_0|+ |(k-1)S+(S\setminus S_0)|\\
&\ge |kS_0|+ k|S+H|-k|H|\\
&\ge k|S_0|-k+1+ k|S+H|-k|H|\ge k|S|-k+1.
 \end{align*}
 In particular, the inequalities used are equalities and hence
 (ii) holds and $|S|=|S+H|-|H|+|S_0|,$ proving (iii). Also, it follows that
 $|kS+H|=|(k-1)S+(S\setminus S_0)|+|H|=k|S+H|-k|H|+|H|,$ proving (iii).
\qed\end{pf}

We can deduce now a Kemperman type result for $kS.$

\begin{cor}

Let $k\ge 2$ be  an integer and let $S$ be a finite  subset of an abelian group $G$ such that $kS$ is aperiodic and $|kS|= k|S|-k+1.$ There is a non-zero subgroup $H$ such $(S\setminus S_0)+H=(S\setminus S_0),$ where $S_0$ is an $H$-component of $S.$  Also, $|kS_0|= k|S_0|-k+1$  and $|k\phi(S)|=k|\phi(S)|-k+1,$  where $\phi :G\mapsto G/H$ denotes
the canonical morphism.
Moreover one of the following holds:
\begin{itemize}
  \item $S_0$ is an arithmetic progression,
\item $k=2$ and $S_0
  =x-((S_0+H)\setminus S_0),$ for some $x.$
\end{itemize}

\label{a-a}
\end{cor}

\begin{pf}
Take a non-zero subgroup $H$ with minimal cardinality
such $k(S+H)=k|S+H|-k|H|+|H|$ and $(S\setminus S_0)+H=(S\setminus S_0),$ where $S_0$ is an $H$-component of $S.$
Notice that $G$ is such a group. Since the period of $kS_0$ is a subgroup of $H$, $kS_0$ is aperiodic and hence $$|kS_0|=k|S_0|-k+1,$$
using the relation $|kS|=k|S|-k+1.$

Observe that $S_0$ can not have a fragment non-zero subgroup
$Q.$ Otherwise we have by Proposition \ref{decomp}, $k(S_0+Q)=k|S_0+Q|-k|Q|+|Q|$ and $(S_0\setminus T_0)+Q=(S_0\setminus T_0),$ where $T_0$ is a $Q$-component of $S_0.$ It would follow that $k(S+Q)=k|S+Q|-k|Q|+|Q|$ and $(S\setminus T_0)+Q=(S\setminus T_0),$ a contradiction. Let $H_0$ be the subgroup generated by $S_0-S_0.$ By Theorem \ref{vosper0}, either (i) holds or one of the following holds:

 \begin{itemize}
   \item $S_0$ is non $2$-separable. We have necessarily $|2S_0| =|H_0|-1.$
Take $a\in S_0$ and put $\{b-a\}=H_0\setminus (2(S_0-a)).$ Necessarily $b-(S_0-a)=H_0\setminus (S_0-a),$ and thus
$b-S_0=H_0+a\setminus (S_0)=(S_0+H_0)\setminus S_0.$

   \item $S_0$ is a $2$-separable Vopser subset.
 We must have $k=2,$ otherwise The condition $|2S_0|\ge \min (|H_0|-1,2|S_0|).$
 But $|H_0|\ge |kS_0|\ge 2|S_0|+|S_0|-1\ge 2|S_0|+1,$ observing that $S_0$ is not an arithmetic progression.
By Kneser's Theorem, $|kS_0|\ge k|S_0|-k+2,$ a contradiction.
Since $|2S_0|=2|S_0|-1$ and since $S_0$ is a Vosper subset, we have necessarily $|2S_0| =|H_0|-1.$
Take $a\in S_0$ and put $\{b-a\}=H_0\setminus (2(S_0-a)).$ Thus, $b-(S_0-a)=H_0\setminus (S_0-a),$ and hence
$$b+a-S_0=(H_0+a)\setminus S_0=(S_0+H_0)\setminus S_0.$$
\end{itemize}
 \qed\end{pf}

 In the above result, the structure of $S$ is completely determined by the structure of $S_0$ and by the structure of $\phi (S).$
 Unfortunately   $k\phi (S)$  is sometimes periodic. In order transform the last result, we investigate the $S,$ where
 $kS$ is periodic and where one element has a unique expressible element. The methods of Kemperman solve very easily this question, as shown in
\cite{hkemp}.

The hyper-atomic approach avoids the last difficulty and lead to a simpler description, as we shall see later.
\section{Universal periods}

Let $T$ and $S$ be finite subsets of an abelian group.
It follows from a result by Balandraud that $|TS|\le |T|+|S|-2$ implies that $T+S$ has a universal period
contained in the atom of $S$. We shall construct a universal period for $kS$ which is bigger in general.

We shall first prove that $S-S$ has a universal period containing the atom if $S$ is not an arithmetic progression and if $|S-S|$ is not very big.

\begin{thm}

Let $r\ge s\ge 1$ be  integers and let $S$ be a finite  subset of an abelian group $G$ and let $H$ be a hyper-atom of $S$.
One of the following holds:
\begin{itemize}
  \item[(i)] $S$ is an arithmetic progression,
  \item[(ii)]  $|sS-rS|\ge \min (|G|-1,(r+s)|S|),$
  \item[(iii)]  The hyper-atom $H$ is a non-zero-subgroup and $sS-rS+H=sS-rS.$
\end{itemize}

\label{a-a}
\end{thm}
\begin{pf}  Assume that (i) and (ii) do not hold. It follows that $S$ is $2$-separable and non-vosperian.
Let $H$ be a hyper-atom of $S.$ By Theorem \ref{vosper0}, $|H|\ge 2.$
By Proposition \ref{decomp}, $S-S+H=S-S.$ Therefore, $sS-rS+H=sS-rS.$  \qed\end{pf}

Proposition \ref{decomp} suggests a very simple  method  giving another universal period for $kS$ containing necessarily Balandraud period.

Let $H$ be a subgroup fragment of $S.$
 An $H$-component $S_0$ of $S$ will be called {\em desertic} component if $|S_0|\le |H|/2.$
 By Proposition \ref{decomp}, the desertic component is unique if it exists. We shall say that
 $S$ is a {\em desert} if it has a desertic component.

%A subset $A_{i+1}$ is an $H_i$-component of $A_{i}$ having a smallest cardinality less than $|H_i|/2$

Given a subset $A$, with $\kappa (A)\le |A|-2.$ We define  a {\em desert} sequence $A_0, \cdots, A_\ell,$
verifying the following conditions:
 \begin{itemize}
   \item  $A_0=A,$
   \item $A_{i+1}$ is a desert  for $0\le i \le \ell-1,$
   \item  $A_\ell$ is not a desert.
 \end{itemize}

Such a sequence exists and is unique, since Proposition \ref{decomp} asserts that $A_i$ is unique
for $1\le i\le \ell.$ The sequence must end since $H_i$ is a finite group with size $<|H_{i-1}|/2.$
The {\em sub-atom} $M$ of $A$ is defied as
 $M=H_\ell$ if $H_\ell$ is non-zero. Otherwise $M=H_{\ell-1}.$ In particular,  the sub-atom is a non-zero subgroup.

\begin{thm}

Let $k$ be a non-negative integer and let $S$ be a finite  subset of an abelian group $G$. If
$|kS|\le k|S|-k,$ then  $$kS+M=kS,$$ where $M$ is the sub-atom of $S.$
\label{subatom} \end{thm}

\begin{pf}
We use the last notations.
The proof is by induction on $\ell$.
 We have $\kappa_1 (S)\le|S|-2
,$ and hence $|H_0|\ge 2.$
By Proposition \ref{decomp}, $(S\setminus S_0)+(k-1)S$ is $H$-periodic. %But $kS\subset kS_0\cup (S\setminus S_0)+(k-1)S$ is $H$-periodic.
 We may assume that $kS_0\cap ((S\setminus S_0)+(k-1)S)=\emptyset,$ otherwise $kS$ is $H_0$-periodic.
Proposition \ref{decomp}, $|kS|=|kS_0|+ |(S\setminus S_0)+(k-1)S|\ge k|S_0|-k+1+ku|H|\ge k|S|-k+1.$
In particular,  $|kS_0|\le k|S_0|-k.$ Notice that $S$ and $S_0$ have the same sub-atom.
By the induction hypothesis $kS_0+M=kS_0$. It follows that $kS+M=kS.$
  \qed\end{pf}

\section{Hyper-atoms and the critical pair Theory}
Applications of hyper-atoms to the critical pair theory where first obtained in \cite{hkemp}.
 A more delicate notion of hyper-atoms was introduced in \cite{hkemp+1}.

\begin{thm}
Let $k\ge 2$ be  an integer and let $S$ be a finite  subset of an abelian group $G$ such that $S$ is not an arithmetic progression, $kS$ is aperiodic and $|kS|= k|S|-k+1.$
 Let $H$ be a hyper-atom of $S$ and let $S_0$ be a smaller $H$-component of $S.$ If $|2S|\neq |G|-1,$ then $|H|\ge 2.$
Moreover, $(S\setminus S_0)+H=(S\setminus S_0)$  and  $|kS_0|=k|S_0|-k+1.$
Also, either $\phi (S)$ is an arithmetic progression or $k=2$ and one of the following holds:
\begin{enumerate}

\item   $S=x-(G\setminus S),$ for some $x.$
\item $(\phi(S)-\phi (S_0))\cap (\phi(S_0)-\phi (S))=\{\phi (0)\}$, where $\phi :G\mapsto G/H$ denotes
the canonical morphism.
\end{enumerate}

\label{a+ahyp}
\end{thm}

\begin{pf}

By Kneser's Theorem and since $2S$ is aperiodic, we have $|2S|=2|S|-1.$
Take an $H$-decomposition $S=S_0\cup \dots \cup S_u.$

Assume first that $S$ is non-$2$-separable. This forces $|2S|=|G|-1.$
Then necessarily $k=2,$ otherwise $3S=G,$ by Lemma \ref{prehistorical}.
Put $2S=G\setminus \{x\}.$ We have clearly $(x-S)\cap S=\emptyset.$ Clearly (1) holds.
Assume now that $S$ is $2$-separable. By Theorem \ref{vosper0}, $|H|\ge 2.$

By Proposition \ref{crdecomp}, $(S\setminus S_0)+H=(S\setminus S_0)$  and  $|kS_0|=k|S_0|-k+1.$

Assume now that $\phi (S)$ is not an arithmetic progression.
By Theorem \ref{hyperatom},
 $\phi (S)$ is  a Vosper subset.

 Thus, $|\phi (G)|-1 < 2|\phi (S)|-1,$ otherwise  $|\phi ((S\setminus S_0)+S)|\ge 2|\phi (S)|-1,$
and hence  $|(S\setminus S_0)+S|\ge 2u|H|+|H|\ge 2|S|,$ a contradiction. Thus, $|\phi (G)| = 2|\phi (S)|-1.$
 In this case, $k=2$ and $2\phi (S)=\phi (G)$. Necessarily, $2\phi(S_0)$
 is uniquely expressible in $2\phi(S).$ In other words $(\phi(S)-\phi (S_0))\cap (\phi(S_0)-\phi (S))=\{\phi (0)\}.$
\qed\end{pf}

\begin{cor}
Let $k\ge 3$ be  an integer and let $S$ be a finite  subset of an abelian group $G$ such that $S$ is not an arithmetic progression, $kS$ is aperiodic and $|kS|= k|S|-k+1.$
 Let $H$ be a hyper-atom of $S$ and let $S_0$ be a smaller $H$-component of $S.$
Then $(S\setminus S_0)+H=(S\setminus S_0)$  and  $|kS_0|=k|S_0|-k+1.$
Moreover  $\phi (S)$ is an arithmetic progression, where $\phi :G\mapsto G/H$ denotes
the canonical morphism.

\label{a+ahyp}
\end{cor}

%%%%%%%%%%%%%%%%%%%%%%%%%%%%%%%%%%%%%%%%%%%%%%%%%%%%%%%%%%%%%%%%%%%%%%%%%%%%%%

%%%%%%%%%%%%%%%%%%%%%%%%%%%%%%%%%%%%%%%%%%%%%%%%%%%%%%%%
%%%%%%%%%%%%%%%%%%%%%%%%%%%%%%%%%%%%%%%%%%%%%%%%%%%%%%%%%%%%%%%%%%%%%%%%%%
%%%%%%%%%%%%%%%%%%%%%%%%%%%%%%%%%%%%%%%%%%%%%%%%%%%%%%%%%%%%%%%%%%%%%%%%%%%%%%%%%
%%%%%%%%%%%%%%%%%%%%%%%
%%%%%%%%%%%%%%%%%%%%%%%%%%%%%%%%%%%%%%%%%%%%%%%%%%%%%%%%%%%%%%%%%%%%%%%%%%%%%%

%%%%%%%%%%%%%%%%%%%%%%%%%%%%%%%%%%%%%%%%%%%%%%%%%%%%%%%%%%%%%%%%%%%%%%%%%%%%%%
%%%%%%%%%%%%%%%%%%%%%%%%%%%%%%%%%%%%%%%%%%%%%%%%%%%%%%%%%%%%%%%%%%%%%%%%%%%%%%

{\bf Acknowledgment}
The author would like to thank M.A. Fiol for stimulating discussions on connectivity 
extra-connectivity

\end{document}